\documentclass[11pt]{article}

\usepackage[latin1]{inputenc}
\usepackage{amsthm,amsmath,amssymb,mathtools}
\usepackage{graphicx}
\usepackage{color}
\usepackage{authblk}
\usepackage{float}
\usepackage[numbers,sort&compress]{natbib}
\usepackage{enumitem}
\usepackage{hyperref}
\hypersetup{colorlinks=true,linkcolor=blue,filecolor=blue,urlcolor=blue,citecolor=cyan,}

\setlist[enumerate]{topsep=0pt,itemsep=-1ex,partopsep=1ex,parsep=1ex}

\usepackage[nameinlink]{cleveref}
\theoremstyle{plain}

\newtheorem{conjecture}{Conjecture}
\newtheorem{theorem}{Theorem}

\newtheorem{proposition}[theorem]{Proposition}

\crefname{corollary}{{Corollary}}{Corollarys}
\crefname{theorem}{Theorem}{Theorems}
\crefname{conjecture}{{Conjecture}}{Conjectures}

\begin{document}

\title{Counterexamples to Clique Immersion Conjecture for Direct Products}
\author{Chuanshu Wu\footnote{Email: wuchuanshu2021@163.com.}~}
 \author{~Zijian Deng\footnote{Email: zj1329205716@163.com (corresponding author).}}


    \affil{Center for Discrete Mathematics\\ Fuzhou University, Fuzhou, 350108, Fujian, China}


\date{\today}
\maketitle

\begin{abstract}
Let \(G\) and \(H\) be graphs, and let \(G\times H\) denote their direct product.
For a graph \(G\), let \(\operatorname{im}(G)\) be the largest integer \(t\)
such that \(G\) contains a \(K_t\)-immersion.  Collins, Heenehan, and
McDonald conjectured that if \(\operatorname{im}(G)=t\) and
\(\operatorname{im}(H)=r\), then
\[\operatorname{im}(G\times H)\ge (t-1)(r-1)+1.\]
We disprove this conjecture by constructing an infinite family of connected
bipartite counterexamples.


\end{abstract}

\section{Introduction}
Throughout this paper, graphs are loopless and simple.
The \emph{degree} of a vertex \(v\) in a graph \(G\), denoted by
\(\deg_G(v)\), is the number of edges incident with \(v\).  
The maximum degree of \(G\), denote by \(\Delta(G)\).
The \emph{chromatic number} \(\chi(G)\) is the minimum number of colours
needed to colour the vertices of \(G\) so that adjacent vertices receive
different colours.

Let \(G\) and \(H\) be graphs.  The \emph{direct product} \(G\times H\)
is the graph with vertex set
$V(G\times H)=V(G)\times V(H).$
Two vertices
$u=(u_1,u_2), v=(v_1,v_2)$
of \(G\times H\) are adjacent if and only if
$u_1v_1\in E(G)$ and $u_2v_2\in E(H).$

Let $G=(V, E)$ be such a graph with vertex set $V(G)$ and edge set $E(G)$.
Let $G,H$ be graphs, $G$ is said to contain a \emph{weak immersion} of $H$ if there exists an injective function $\phi: V(H)\rightarrow V(G)$, satisfying the following:	

\begin{itemize}
        \item for every edge $uv\in E(H)$, there is a path in $G$, denoted $P_{uv}$, with endpoints $\phi(u)$ and $\phi(v)$.
        \item the paths in $\{P_{uv} | uv\in E(H)\}$ are pairwise edge disjoint.
\end{itemize}

The vertices $\{\phi(v):v\in V(H)\}$ are called the \emph{terminal vertices} (also known as \emph{branch vertices}) of the immersion, and call internal vertices of the path $P_{uv}$ (where $u,v \in V(G)$) the \emph{pegs} of the $H$-immersion.
If terminal vertices do not appear as interior vertices on path $P_{uv}$, we say that $G$ contains a \emph{strong immersion} of $H$.
This paper is only concerned with weak immersion, and from now  on we will omit ``weak'', and simply say ``immersion''.
The \emph{immersion number} of a graph $G$ be denoted as $im(G)$, is the largest value $t$ for which $G$ has a $K_{t}$-immersion.
$K_{t}$ is a \emph{complete graph} on $t$ vertices in which any two vertices are adjacent.

The immersion relation has received a considerable attention in recent years, especially after Robertson and Seymour \cite{RS2010} proved that it is a well-quasi-order, just like the minor relation.
Much of this attention has been directed toward the conjecture of Lescure and Meyniel \cite{LM1989}, which serves as an analogue to Hadwiger's conjecture \cite{H1943} for strong immersion.
It states that every graph $G$ with a chromatic number $\chi(G)\geq t$ contains the complete graph $K_{t}$ as a strong immersion.
Additionally, Abu-Khzam and Langston \cite{AL2003} proposed a conjecture for weak immersion.

Haj\'os' conjectured that every graph $G$ with a chromatic number $\chi(G)\geq t$ contains the complete graph $K_{t}$ as a subdivision.
The counterexamples that Catlin \cite{C1979} used to disprove Haj\'os' conjecture are the lexicographic product of cliques and cycles.

Collins et al. \cite{CHM2023} studied clique immersions in the four standard
graph products: the lexicographic product \(G\circ H\), the Cartesian
product \(G\square H\), the direct product \(G\times H\), and the strong
product \(G\boxtimes H\).  They first obtained general upper bounds on
the immersion number of \(G\ast H\) in terms of the maximum degree of the
product graph.  For lower bounds, they considered the following natural
question: if $\operatorname{im}(G)=t $ and $\operatorname{im}(H)=r,$
is it true that $\operatorname{im}(G\ast H)\ge \operatorname{im}(K_t\ast K_r)?$
In~\cite{CHM2023}, the corresponding lower bounds for the lexicographic
product and the Cartesian product were completely determined.  However,
the direct product and the strong product have substantially different
structures, which makes the corresponding lower-bound problems more
difficult. 
For the direct product, Collins et al. proposed the following conjecture. 

\begin{conjecture}(Collins, Heenehan, and McDonald \cite{CHM2023})\label{Con 2}
Let $G$ and $H$ be graphs with $im(G)=t$ and $im(H)=r$, then $im(G \times H) \geq (t-1)(r-1)+1$.
\end{conjecture}

They proved that \cref{Con 2} holds in several special cases, including
the case where \(H\) contains \(K_r\) as a subgraph, and the case where the
paths corresponding to the clique immersions in \(G\) and \(H\) all have
the same parity.
For the strong product, Collins et al.~\cite{CHM2023} conjectured that
$\operatorname{im}(G\boxtimes H)\ge tr$
whenever \(\operatorname{im}(G)=t\) and \(\operatorname{im}(H)=r\).  Wu
and Deng~\cite{WD2025} verified this conjecture for several classes of
graphs admitting suitable clique immersions.  In particular, they proved
the conjecture when \(H\) contains \(K_r\) as a subgraph.

In this paper, we construct a family of counterexamples disproving \cref{Con 2}.

\section{ The Counterexample}

We define a class of graphs \(B_t^{(p,q)}\), where \(t\ge 3\),
\(p+q=t\), and \(1\le p\le t-1\).

Let
\[
A=\{a_1,\ldots,a_p\},\qquad B=\{b_1,\ldots,b_q\}
\]
be a partition of a set of \(t\) branch vertices.  The graph
\(B_t^{(p,q)}\) is obtained from the complete graph on \(A\cup B\) by
subdividing exactly once every edge whose endpoints both lie in \(A\) or
both lie in \(B\), while leaving every edge between \(A\) and \(B\)
unchanged.

This construction has several immediate properties:
\[|V(B_t^{(p,q)})|=t+\binom{p}{2}+\binom{q}{2},\quad |E(B_t^{(p,q)})|=pq+2\binom{p}{2}+2\binom{q}{2},\]
and
\[\Delta(B_t^{(p,q)})=t-1.\]

More importantly, every branch vertex has degree \(t-1\), whereas every
newly introduced subdivision vertex has degree \(2\).  Since each edge
with both endpoints in the same part is replaced by a path of even length
\(2\), while each edge between \(A\) and \(B\) is kept as a path of odd
length \(1\), the graph \(B_t^{(p,q)}\) is bipartite.

To describe its bipartition explicitly, let \(S_A\) denote the set of
subdivision vertices introduced on edges with both endpoints in \(A\), and
let \(S_B\) denote the set of subdivision vertices introduced on edges with
both endpoints in \(B\).  Then \(B_t^{(p,q)}\) has bipartition
\[(A',B')=(A\cup S_B,\; B\cup S_A).\]
Indeed, every edge of \(B_t^{(p,q)}\) has one endpoint in \(A'\) and the
other in \(B'\).  Moreover,
\[|A'|=p+\binom{q}{2},\qquad|B'|=q+\binom{p}{2}.\]
Since the original complete graph on \(A\cup B\) is connected, and
subdividing edges does not destroy connectivity, the graph
\(B_t^{(p,q)}\) is connected.  Hence \(B_t^{(p,q)}\) is a connected
bipartite graph.

Moreover, there is an explicit \(K_t\)-immersion in \(B_t^{(p,q)}\).
Take the branch vertices \(A\cup B\) as the set of terminals.  For each
pair of terminals lying in the same part, use the corresponding path of
length \(2\); for each pair of terminals lying in different parts, use the
corresponding edge.  By construction, these paths are pairwise
edge-disjoint.  Hence
\[\operatorname{im}(B_t^{(p,q)})\ge t.\]

Conversely, since
\[\Delta(B_t^{(p,q)})=t-1,\]
the graph \(B_t^{(p,q)}\) cannot contain a \(K_{t+1}\)-immersion.  Indeed,
in any \(K_{t+1}\)-immersion, each terminal must be incident with at least
\(t\) distinct edges, one for each of the edge-disjoint paths from that
terminal to the other \(t\) terminals.  This is impossible in a graph of
maximum degree \(t-1\).  Therefore
\[\operatorname{im}(B_t^{(p,q)})\le t.\]
Combining the two inequalities, we obtain
\[\operatorname{im}(B_t^{(p,q)})=t.\]

Now let
\[G=B_t^{(p,q)},\qquad H=B_r^{(s,u)},\]
where \[p+q=t,\qquad s+u=r,\]
and assume that \[p,q,s,u\ge 2.\]
Set \[k=(t-1)(r-1)+1.\]
This is precisely the lower bound predicted by the direct product
conjecture for \(\operatorname{im}(G\times H)\).

We first record a few key observations.  Since both \(G\) and \(H\) are
connected bipartite graphs, their direct product \(G\times H\) has exactly
two connected components.  Moreover, degrees in the direct product satisfy
\[\deg_{G\times H}(x,y)=\deg_G(x)\deg_H(y).\]
Since \(p,q,s,u\ge2\), we have \(t,r\ge4\).  Hence the only vertices of
\(G\times H\) attaining degree \[(t-1)(r-1)=k-1\]
are those of the form \((x,y)\), where both \(x\) and \(y\) are branch
vertices.  Indeed, every non-branch vertex in either factor has degree
\(2\), and therefore any product vertex with at least one non-branch
coordinate has degree strictly less than \((t-1)(r-1)\).

Now every terminal of a \(K_k\)-immersion must have degree at least
\(k-1\).  Consequently, if \(G\times H\) contained a \(K_k\)-immersion,
then all its \(k\) terminals would have to be branch-by-branch vertices.
Furthermore, since all terminals in a clique immersion must lie in the
same connected component, these \(k\) branch-by-branch vertices would all
have to lie in a single connected component of \(G\times H\).

It remains to count how many branch-by-branch vertices occur in each
component.  Let \(A\cup B\) be the branch bipartition of \(G\), and let
\(C\cup D\) be the branch bipartition of \(H\), where
\[|A|=p,\qquad |B|=q,\qquad |C|=s,\qquad |D|=u.\]
The two connected components of \(G\times H\) contain, respectively,
\[|A||C|+|B||D|=ps+qu\] and \[|A||D|+|B||C|=pu+qs\]
branch-by-branch vertices.

We now compare these quantities with \(k\).  Since
\[k=(t-1)(r-1)+1,\]
we have
\[k-(ps+qu)=(t-1)(r-1)+1-(ps+qu)=(p-1)(u-1)+(q-1)(s-1),\]
and 
\[k-(pu+qs)=(t-1)(r-1)+1-(pu+qs)=(p-1)(s-1)+(q-1)(u-1).\]
Because \(p,q,s,u\ge2\), both right-hand sides are at least \(2\).  Hence
\[ps+qu<k \qquad\text{and}\qquad pu+qs<k.\]
Equivalently,
\[\max\{ps+qu,\;pu+qs\}<k.\]
Thus neither connected component of \(G\times H\) contains enough
branch-by-branch vertices to serve as the terminal set of a
\(K_k\)-immersion.  Therefore
\[\operatorname{im}(G\times H)<k.\]

We summarize this in the following proposition.

\begin{proposition}
Let $G=B_t^{(p,q)}$ and $H=B_r^{(s,u)},$
where $p+q=t,s+u=r,$ and \(p,q,s,u\ge2\).  Then
$\operatorname{im}(G)=t$ and $\operatorname{im}(H)=r,$
but
\[
\operatorname{im}(G\times H)<(t-1)(r-1)+1.
\]
\end{proposition}

Consequently, \cref{Con 2} fails for
\[
B_t^{(p,q)}\times B_r^{(s,u)}.
\]


\section*{Declaration}
\noindent$\textbf{Conflict~of~interest}$
	The authors declare that they have no known competing financial interests or personal relationships that could have appeared to influence the work reported in this paper.
	
	\noindent$\textbf{Data~availability}$
	Data sharing not applicable to this paper as no datasets were generated or analysed during the current study.


\begin{thebibliography}{100}

\baselineskip12pt

\bibitem{AL2003} F. N. Abu-Khzam, M. A. Langston,
Graph coloring and the immersion order,
\textit{Computing and Combinatorics, in: Lecture Notes in Computer Science} \textbf{2697} (2003) 394--403.

\bibitem{C1979} P. A. Catlin,
Haj\'os' graph-coloring conjecture: variations and counterexamples,
   \textit{J. Combin. Theory, Ser. B} \textbf{26} (1979) 268--274.
   
\bibitem{CHM2023} K. L. Collins, M. E. Heenehan, J. McDonald,
Clique immersion in graph products,
\textit{Discrete Math.} \textbf{346} (2023) 113421.

\bibitem{H1943} H. Hadwiger,
\"{u}ber eine Klassifikation der Streckenkomplexe,
\textit{Vierteljschr. Naturforsch. Ges. Z\"{u}rich} \textbf{88} (1943) 133--142.

\bibitem{LM1989} F. Lescure, H. Meyniel,
On a problem upon configurations contained in graphs with given chromatic number,
\textit{Ann. Discrete Math.} \textbf{41} (1989) 325--331.

\bibitem{RS2010} N. Robertson, P. Seymour,
Graph minors XXIII, Nash-Williams' immersion conjecture,
\textit{J. Combin. Theory, Ser. B} \textbf{100} (2010) 181--205.

\bibitem{WD2025} C. Wu, Z. Deng,
A note on clique immersion of strong product graphs,
\textit{Discrete Math.} \textbf{348} (2025) 114237.

\end{thebibliography}
\end{document}